\documentclass{article}

\usepackage{arxiv}

\usepackage[utf8]{inputenc} 
\usepackage[T1]{fontenc}    
\usepackage{hyperref}       
\usepackage{url}            
\usepackage{booktabs}       
\usepackage{amsfonts}       
\usepackage{nicefrac}       
\usepackage{microtype}      
\usepackage{graphicx}
\usepackage{doi}
\usepackage{tabularx} 

\newcommand\myatop[2]{\genfrac{}{}{0pt}{}{#1}{#2}}

\usepackage{caption}

\usepackage{blindtext}

\usepackage[leqno]{amsmath}
\usepackage{amsthm}
\usepackage{amssymb}
\usepackage{booktabs} 
\usepackage{siunitx}
\usepackage{physics}
\usepackage{multirow}

\usepackage{todonotes}
\def\bX{\mathbf{X}}

\def\bA{\mathbf{A}}
\def\bB{\mathbf{B}}
\def\bF{\mathbf{F}}

\def\bQ{\mathbf{Q}}
\def\bS{\mathbf{S}}
\def\bV{\mathbf{V}}
\def\bR{\mathbf{R}}

\def\bTheta{\mathbf{\Theta}}
\def\btheta{\boldsymbol{\theta}}
\def\onevec{{1}}

\newtheorem{theorem}{Theorem}
\newtheorem{lemma}{Lemma}
\newtheorem{definition}{Definition}

\usepackage[backend=bibtex,style=authoryear,natbib]{biblatex}
\usepackage{csquotes}
\title{Efficient semidefinite bounds for multi-label discrete graphical models}

\author{

    \href{https://orcid.org/0000-0001-7883-4185}{\includegraphics[scale=0.06]{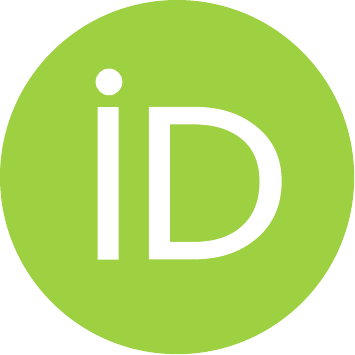}\hspace{1mm}Valentin Durante}\\
	Université Fédérale de Toulouse\\
	ANITI, INRAE, UR 875\\
	Toulouse, France\\
	\texttt{valentin.durante@inrae.fr}\\
	
	\And
	
	\href{ https://orcid.org/0000-0002-3727-6698}{\includegraphics[scale=0.06]{orcid.pdf}\hspace{1mm}George Katsirelos}\\
	UMR MIA-Paris, INRAE, AgroParisTech\\
	Univ. Paris-Saclay\\
	75005 Paris, France\\
	\texttt{gkatsi@gmail.com}\\
	
    \And
    
    \href{https://orcid.org/0000-0001-6049-3415}{\includegraphics[scale=0.06]{orcid.pdf}\hspace{1mm}Thomas Schiex}\\
	Université Fédérale de Toulouse\\
	ANITI, INRAE, UR 875\\
	Toulouse, France\\
	\texttt{thomas.schiex@inrae.fr}\\
}

\renewcommand{\headeright}{}
\renewcommand{\undertitle}{}

\hypersetup{
pdftitle={Efficient semidefinite bounds for multi-label discrete graphical models},
pdfsubject={q-bio.NC, q-bio.QM},
pdfauthor={Valentin Durante, George Katsirelos, Thomas Schiex},
pdfkeywords={Weighted Constraint Satisfaction, Convex Optimization, Low rank, Graphical Models, Discrete Optimization}
}

\begin{document}
\maketitle

\begin{abstract}
 By concisely representing a joint function of many variables as the combination of small functions, discrete graphical models (GMs) provide a powerful framework to analyze stochastic and deterministic systems of interacting variables. One of the main queries on such models is to identify the extremum of this joint function. This is known as the Weighted Constraint Satisfaction Problem (WCSP) on deterministic Cost Function Networks and as Maximum a Posteriori (MAP) inference on stochastic Markov Random Fields. 
  
  Algorithms  for approximate WCSP inference typically rely on local consistency algorithms or belief propagation. These methods are intimately related to linear programming (LP) relaxations and often coupled with reparametrizations defined by the dual solution of the associated LP.
  
  Since the seminal work of Goemans and Williamson, it is well understood that convex SDP relaxations can provide superior guarantees to LP. But the inherent computational cost of interior point methods has limited their application. The situation has improved  with the introduction of non-convex Burer-Monteiro style methods which are well suited to handle the SDP relaxation of combinatorial problems with binary variables (such as MAXCUT, MaxSAT or MAP/Ising).
  
  We compute low rank SDP upper and lower bounds for discrete pairwise graphical models with arbitrary number of values and arbitrary binary cost functions by extending a Burer-Monteiro style method based on row-by-row updates. We consider a traditional dualized constraint approach and a dedicated Block Coordinate Descent approach which avoids introducing large penalty coefficients to  the formulation. On increasingly hard and dense WCSP/CFN instances, we observe that the BCD approach can outperform the dualized approach and provide tighter bounds than local consistencies/convergent message passing approaches.
\end{abstract}


\section{Introduction, Definition} 

Graphical models~\citep{cooper2020graphical} are descriptions of decomposable multivariate functions. A variety of frameworks in Computer Science, Logic, Constraint Satisfaction/Programming, Machine Learning, Statistical Physics, Operations Research and Artificial Intelligence define specific sub-classes of graphical models.

\newpage
\begin{definition}{(Graphical Model)}
A Graphical Model $M=\langle X,\Phi\rangle$ with co-domain B and a combination operator $\oplus$ is defined by :
\begin{itemize}
    \item a sequence of n variables $X$.
    \item a set of functions $\Phi$. Each function $\theta_S \in \Phi$ is a function from $D_S \rightarrow B$. S is called the scope of the function and $|S|$ its arity. $D_S$ is the cartesian product of the domain of the variables within the scope. 
\end{itemize}
$M$ defines a joint function:
\begin{center}
$\Theta_M(v) = \displaystyle \bigoplus_{\theta_S \in \Phi} \theta_S(v[S])$
\end{center}
Where $v[S]$ is the projection of $v$ on $S = scope(\theta_S)$: the sequence of values that the variables in $S$ takes in $v$. 
\end{definition}
We consider additive numerical GMs, aka energy-based models that are used to describe infinite-valued or bounded cost functions in Constraint Programming, using Valued Constraint and Cost function Networks~\citep{schiex1995valued,cooper2010soft} or to describe probability distributions as Markov Random Fields or Bayesian Networks in the stochastic case~\citep{koller2009probabilistic}. GMs have many areas of applications, in image processing, computational biology, resource management, configuration, etc.

For our work we are interested in solving optimization queries. A GM $M=\langle X,\Phi\rangle$ is specified by a set of variables $X=\{x_1,\ldots,x_n\}$, where $x_i$ can take $d_i$ possible values and a set $\Phi$ of potential or cost functions $\theta_{S}$ involving variables in $S\subset X$, the scope of $\theta_S$. We note $d=\sum_{i=1}^n d_i$ the total number of values. An additive graphical model defines a joint function $\Theta_M$ over $X$ as $\Theta_M=\sum_{\theta_{S}\in \Phi} \theta_{\bS}$. Deterministic GMs use non negative integer functions and a saturating arithmetic. Stochastic models use real $\theta_{S}$ functions, exponentiation and normalization to define a probability distribution.  When all scopes have  cardinality less than or equal to 2, the model is pairwise and its graph has vertices in $X$ and edges defined by scopes.

Our optimization query is to find an assignment of the discrete variables $v \in D_X$ that minimizes the joint function.
\begin{equation*}
    \min_{v\in D_X} \displaystyle \sum_{\theta_S \in \Theta} \theta_S(v[S])
\end{equation*}
This problem is known as the Weighted Constraint Satisfaction Problem (WCSP). In the case of a pairwise GM, finding an assignment of the variables that minimizes the joint function can be expressed as an Integer Linear Program (ILP).
\begin{equation}
    \begin{aligned}
      \mbox{Minimize} & \sum_{i,a} \theta_i(a)\cdot x_{ia} + &\sum_{\myatop{\theta_{ij}\in \Phi}{a\in D_i, b\in D_j}} \theta_{ij}(a,b)\cdot y_{iajb} \ \ \mbox{ such that }\nonumber  \\[3mm]
      &\sum_{a\in D_i} x_{ia} = 1  &\forall i \in \{1,\ldots,n\} \nonumber\\
      & \sum_{b\in D_j} y_{iajb} = x_{ia}  &\forall \theta_{ij}\in \Phi, \forall a \in D_i\nonumber\\
      &\sum_{a\in D_i} y_{iajb} = x_{jb}  &\forall \theta_{ij}\in \Phi, \forall b \in D_j\nonumber\\
      & x_{ia}\in\{0,1\}  &\forall i\in \{1,\ldots,n\}\nonumber
    \end{aligned}
\end{equation}
Here, $x_{ia}$ is a binary variable which is 1 if variable $x_i \in X$ takes value $a \in D_i$, 0 otherwise. The second line describes the "exactly-one" constraints, a GM variable can take exactly one value within its domain. $y_{iajb}$ are 0/1 variables that encode assignments of pairs of variables, for all the binary cost functions. The second and third constraints ensure that those pair variables are assigned in a way which is consistent with the assignments of the $x_{ia}$ variables.

Except for specific polynomial classes that include tree-width bounded models or models defined by submodular functions~\citep{cooper2020graphical}, the problem of finding an  assignment of the variables in $X$ that optimizes the joint function $\Theta_M$ is decision NP-complete. While exact Branch and Bound based algorithms exist for the problem,  a vast array of polynomial heuristics have been developed for the most challenging instances, that remain out of reach of exact algorithms. Relatively small very hard instances can be easily produced as random problems. Such hard instances also easily pop-up when Machine Learning algorithms are used to estimate GMs from data using regularized loglikelihood approaches, when tuning the regularizing penalty~\citep{brouard2020pushing}. 

While the integer linear program is intractable in general, one can derive a Linear Programming (LP) relaxation of the original problem. The decision variables $x_{ia}$ will now take continuous value in the interval $[0,1]$. This LP relaxation, known as the local polytope~\citep{Wer07}, leads to a simpler polytope to optimize over, with a tractable number of constraints. If the LP relaxation returns an integer valued solution, this is an optimal solution of the original ILP~\citep{wainwright2005map}. However, even though solving LP is tractable, the number of constraints makes direct LP solvers inefficient on large instances. Efficient heuristics were developed to address this inefficiency.

A first family of heuristics derives from the seminal convergent arc consistency~\citep{waltz1975understanding} (or Boolean Constraint Propagation, BCP) algorithm introduced for Boolean graphical models (Constraint Networks). This algorithm was generalized by Pearl to numerical GMs as Loopy Belief Propagation~\citep{pearl1988probabilistic}, at the cost of losing convergence. Convergence was restored in Soft  Arc Consistency algorithms~\citep{schiex2000arc,cooper2010soft} for the WCSP/CFN, culminating in Virtual Arc Consistency (VAC)~\citep{cooper2008virtual}. Independetly to the development of VAC, convergent belief propagation algorithms such as \mbox{TRW-S}~\citep{kolmogorov2005convergent} were introduced, which witness dual bounds by a potentially suboptimal LP dual solution, which can be interpreted as a reparametrizations of the MAP/MRF. Fixpoints of algorithms such as TRW-S share the same properties as those of VAC, but VAC
additionally applies BCP on infinite-valued energies. These algorithms offer lower bounds for minimization. They also heuristically produce primal (discrete) heuristic solutions. Taken together, these offer post-hoc guarantees through an optimality gap, making such algorithms useful in time-constrained scenarios.

In many cases, these bounds correspond to high quality solutions and a small optimality gap. On the hardest instances, however, they produce primal solutions of high cost and poor lower bounds, failing to provide reasonable guarantees: efficient tighter bounds are needed. One natural direction to obtain such bounds is to use higher levels of the Sherali-Adams hierarchy \citep{SA90}, which in practical terms means extending Arc Consistency and BP to a larger set of variables. This idea is captured by $k$-consistency for Boolean GMs~\citep{journals/cacm/Freuder78,journals/jacm/Freuder82} and CFNs~\citep{nguyen2017triangle}, or GBP~\citep{yedidia2000generalized} and convergent variants~\citep{sontag2012tightening} for MAP/MRF. The associated time and space costs quickly become extreme however, as moving up one step in the hierarchy multiplies time and space complexity by a factor equal to the average domain size. A different approach is to use a different hierarchy, namely to try to exploit SDP bounds which provide superior guarantees~\citep{goemans1995improved} compared to approximate LP bounds.

However, SDP has found limited applications in solving optimization queries in graphical models. The main reason lies in the fact that Goemans and Williamson's SDP relaxation of MAXCUT, a specific case of MAP/MRF with binary variables, requires $\theta(n^2)$ variables and memory resulting in an $O(n^6)$ complexity for interior point methods. The algorithm doesn't scale well and becomes impractical except on small very hard problems. It very quickly trades off too much in speed in favor of strength of inference. For this reason approximate LP methods such as those we describe above are preferred because they achieve a better trade off.

This has changed with the non-convex low-rank approach proposed by \citep{burer2003nonlinear}. They exploited the result of \citep{barvinok1995problems} and \citep{pataki1998rank}, who showed that the rank of solutions of SDP problems is bounded by $O(\sqrt{n})$. Burer and Monteiro use a low-rank factorization $\bX = \bV\bV^T$ of the semidefinite matrix $\bX$ in the SDP relaxation of these combinatorial problems. It has since been observed that, in practice, the rank of solutions is even lower, some software using a constant $O(1)$ rank. Then, the number of variables and the memory requirements can be reduced to $O(nr)$ where $r$ is the decomposition rank used. As a side benefit, the decomposition guarantees the semidefiniteness of the matrix $\bX$. This comes at the cost of addressing a non-convex optimization problem.

Combined with row-by-row updates~\citep{javanmard2016phase,wang2017mixing}, Riemannian gradient descent~\citep{javanmard2016phase,mei2017solving} or Riemannian trust-region methods~\citep{absil2007trust}, the Burer-Monteiro approach provides empirically very efficient methods to solve SDP relaxations of combinatorial optimization problems such as MAXCUT, MAXSAT~\citep{wang2019low} or MAP over Ising or Potts models~\citep{pabbaraju2020efficient}. The row-by-row update methods, with their efficient $O(nr)$ updates are especially attractive. 

Potts model MRFs, where pairwise functions are weighted identity matrices can still be directly dealt by using simple row-by-row updates~\citep{pabbaraju2020efficient}. In contrast, arbitrary pairwise GMs, with non binary variables and arbitrary functions need to explicitly deal with one additional constraint for each original variable of the considered GM. In the usual one-hot encoding of finite domain variables, a GM variable $x_k$ with $d_k$ values is encoded as $d_k$ binary variables, each representing whether a value is used or not. The extra constraint specifies that a variable is assigned exactly one value. We add these constraints to the SDP relaxation and dualize them with a  penalty of sufficient magnitude to guarantee they will be enforced~\citep{lasserre2016max}. Finally, we may try to get rid of these constraints altogether by shifting from row-by-row updates to Block Coordinate (BC) updates on blocks that include all the $rd_k$ variables that represent the $d_k$ binary variables in the low-rank matrix $\bV$. This approach requires however to have efficient BC updates.

We show that such BC updates can be performed efficiently while preserving convergence. We implement the two approaches described above and compare them with \mbox{TRW-S} and VAC.

\section{SDP relaxation for pairwise GMs}
\label{headings}

Matrices are denoted with a capital boldface font and all vectors are column vectors. For matrix $\bA$, $A_{ij}$ represents the entry at the $i$-th row and $j$-th column of $\bA$. $A_i$ represents its $i$-th row as a column vector. For a vector $g$, $||g||$ denotes its Euclidean norm. The vector defined by the diagonal of $\bA$ is $\textrm{diag}(\bA)$, the Frobenius inner product of matrices $\bA,\bB$ is denoted as $\langle \bA,\bB\rangle$. For a vector g, $\textrm{diag}(g)$ denotes the diagonal matrix with diagonal $g$.  

A pairwise GM $M=\langle X,\Phi\rangle$ with associated graph $(X,E)$ defines the joint function $\Theta_M = \sum_{(x_i,x_j) \in E} \theta_{ij} +\sum_{x_i\in X} \theta_i + \theta_\varnothing$. We assume that the functions $\theta_S\in \Theta$ are described as tensors, here matrices, vectors and a scalar respectively. For optimization, the constant term $\theta_\varnothing$ can be ignored. Before building an SDP relaxation, the problem is reduced to a quadratic form by representing every finite domain variable $x_k\in X$ as a vector $b_k$ of $d_k$ Boolean variables. The $j^{th}$ element of $b_k$ will take value $1$ when variable $x_k = j$ and value $0$ otherwise. The value of $\theta_{ij}(x_i,x_j)$ is $b_i^T \bTheta_{ij} b_j$ and the value of $\theta_i(x_i)$ is then $b_i^T \theta_i$. If we denote by $b$ the stacked vector $b^T = [b_i^T,\ldots,b_n^T]$, we then have $\Theta_M(x) = b^T \bTheta^{\mathbb{B}} b + b^T \theta^{\mathbb{B}}$ where $\theta^{\mathbb{B}}$ is the stacked vector of all $\theta_i$ and $2\times \bTheta^{\mathbb{B}}$ is a symmetric block matrix having block $\bTheta_{ij}$ when $(x_i,x_j)\in E$ and the $d_i\times d_j$ zero matrix otherwise. Note that $\bTheta^{\mathbb{B}}$ is sparse by block, according to the connectivity of the graph of $M$. More importantly, $\bTheta^{\mathbb{B}}$ has zero blocks on its diagonal since no pairwise function in $\Phi$ connects a variable $x_k$ to itself.

When optimizing over $b$ instead of $x$, one must enforce the fact that only one element of every Boolean vector $b_k$ is set to $1$, \emph{i.e.}, that $\forall x_k\in X, \sum b_k =1$. This set of constraints can be gathered in a linear constraint $\bA b = \onevec_n$ where $\bA$ is a Boolean matrix in $\mathbb{B}^{n\times d}$ with $A_i^T = [0_{d_1}^T \; \cdots \; 0_{d_{i-1}}^T \; 1_{d_i}^T \;  0_{d_{i+1}}^T \; \cdots \; 0_{d_n}]$ and where $\onevec_n$ is an $n$-vector of $1$. Finally, minimizing $\Theta_M$ becomes equivalent to solving the following constrained quadratic program:

\begin{equation}
   \min_{b \in \mathbb{B}^d}  \; \; b^T \bTheta^{\mathbb{B}} b + b^T \btheta^{\mathbb{B}} \; \; \textrm{ s.t } \; \; \bA b = \onevec_n\\
    \label{eq2}
\end{equation}

To build an SDP relaxation of the quadratic program above, centered $\{-1,1\}$ variables are used and relaxed to norm-1 vectors of dimension $n$. The problem above can be reduced to such formulation by using the $\{-1,1\}$ variables $c= 2b - 1_d$. After this variable change, denoting $\bTheta = \bTheta^{\mathbb{B}}/4$, $\theta = (\theta^{\mathbb{B}} + \bTheta^{\mathbb{B}} 1_d)/2$,  $\bF = \bA/2$ and $u = 1_n - \bF 1_d/2$, the problem becomes:

\begin{equation}
   \min_{c \in \{-1,1\}^d}  \; \; c^T \bTheta c + c^T \theta \; \; \textrm{ s.t } \; \; \bF c = u
    \label{eq9}
\end{equation}
Here, the initial exactly-one constraints $\sum b_i =1$ become $2F_i c = 2 - d_i$, with $F_i \in \mathbb{R}^d$, $2F_i^T = [0_{d_1}^T \; \cdots \; 0_{d_{i-1}}^T \; 1_{d_i}^T \;  0_{d_{i+1}}^T \; \cdots \; 0_{d_n}]$. This change of variables creates a constant term in the objective value which can be ignored for optimization purposes. In this formulation, the matrix $\bTheta$ has the exact same sparsity as $\bTheta^{\mathbb{B}}$.

\subsection{Lassere dualization of linear constraints}

To get rid of the linear constraints in (\ref{eq9}) and get a quadratic program, a usual approach is to dualize the linear constraint with a penalty $\rho$ that must be properly chosen~\citep{lasserre2016max}. This reduces the problem to a pure MAXCUT problem on which the Burer-Monteiro approach can be directly applied using its efficient row-by-row updates~\citep{wang2017mixing}. The relaxed problem then becomes:

\begin{equation}
    \displaystyle \min_{c \in \{-1,1\}^d} \; \; c^T \bTheta c + c^T \theta  + (2\rho + 1)||Fc - u||^2 \label{eq4}
\end{equation}

As soon as $\rho \geqslant \max\{|c^T \bTheta c + c^T \theta| \; : \; c \in \{-1,1\}^d \}$, the solution of (\ref{eq4}) and (\ref{eq2}) are the same~\citep{lasserre2016max}. One can note at this point that if any potential function in the GM at hand contains functions of large amplitude, the above property requires large values of $\rho$. 

In order to convert linear terms to quadratic terms (for homogenization), we introduce the extended vector $e^{T} = [c~1] \in \{-1,1\}^{d+1}$. Then, (\ref{eq4}) equivalently asks to minimize $e^T \bQ e$ where $e \in \{-1,1\}^{d+1}$ and $\bQ$ (symmetric) is:
\[
  \begin{pmatrix}
   (2\rho + 1)F^T F + \bTheta  & \frac{1}{2} (\theta^T -2(2\rho+1)u^TF)^T \\
    & \\
    \frac{1}{2} (\theta^T -2(2\rho+1)u^T F)\!\!\!\! & (2 \rho + 1)u^Tu
   \end{pmatrix} \]

The usual SDP relaxation of the MAX-CUT problem can be written using the rank 1 matrix variable $\bX = e e^T \in \mathbb{R}^{(d+1)\times (d+1)}$. Dropping the rank-1 constraint, the relaxation is:
\begin{equation}
    \displaystyle \min_{\bX}\{ \langle \bQ,\bX\rangle : \bX \succeq 0; \; X_{ii} = 1, \; i = 1,\ldots,d+1 \}
    \label{eq7}
\end{equation}

\citep{burer2003nonlinear} introduced the idea of using a low-rank factorization of $\bX$ to solve the SDP (\ref{eq7}). This factorization was motivated by a proof by~\citep{barvinok1995problems} and~\citep{pataki1998rank} of the following theorem:
\begin{theorem}
There exists an optimum $\bX^*$ of (\ref{eq7}) with rank $r$ such that $\frac{r(r+1)}{2} \leqslant d + 1$.
\end{theorem}

Every positive semidefinite matrix of rank $r$ can be factorized as a product of two rank $r$ matrices: $\bX = \bV \bV^T, \; \bV \in \mathbb{R}^{(d+1) \times r}$. Then, (\ref{eq7}) becomes:
\begin{equation}
    \min_{\bV \in \mathbb{R}^{(d+1) \times r}} \big\langle \bQ,\bV \bV^T \big\rangle \text{ s.t } ||V_i|| = 1 \; , \; i = 1,\ldots,d+1
    \label{eq8}
\end{equation}
With $V_i \in \mathbb{R}^r$ the row vector $i$ of $\bV$, the number of variables shifts from $(d+1)^2$ to $r(d+1)$ and the constraint $\bX \succeq 0$ becomes implicit as $\bV \bV^T$ is always positive semidefinite. Several approaches exploit this idea. We use the mixing method~\citep{wang2017mixing}. It does efficient cyclic $O(rd)$ updates of the row vectors $V_i$, all other row vectors being fixed:
$$V_i = -\frac{g_i}{||g_i||}, \textrm{ where } g_i = \displaystyle \sum_{j=1}^{d+1} Q_{i,j} V_j \quad \forall i = 1,\ldots,d+1$$

The mixing method is known to recover the optimum of the convex SDP (\ref{eq7}) as long as $r > \sqrt{2*(d+1)}$. 
One issue of the mixing method in the presence of dualized constraints generated by Lassere's approach is that the row-by-row updates will be sensitive to the magnitude of the penalty $\rho$. As we observed previously, $\rho$ may be large if the input GM has large terms. The large value of $\rho$ will create large coefficients in $\bQ$ and then the value of $g_i$ will be dominated by few terms, possibly slowing down convergence. It is important to notice that if the solution to (\ref{eq7}) nullifies the quadratic penalty term induced by the Lasserre dualization, it will implicitly satisfy the so-called \textit{gangster} constraint~\citep{zhao1998semidefinite}. If we consider the Boolean vector $b_k$ that corresponds to the variable $x_k \in X$, the gangster constraint ensures that the off-diagonal entries of the matrix $b_kb_k^T$ are all zeros. Indeed, $b_{k}^i*b_{k}^j = 0  \; \forall \; i \ne j$, since at least one of the two terms is equal to zero. This can tighten the SDP relaxation bound.

\subsection{SDP bound using block coordinate descent}

We now consider an approach that exploits the specific structure of the matrix $\bQ$ in the constrained quadratic formulation (\ref{eq9}). Using the extended vector $e$,  the rank-1 matrix variable $\bX= e e^T$, dropping the rank-1 constraint and using a low rank factorization $\bX = \bV\bV^T$ as above, the quadratic optimization problem (\ref{eq9}) can be relaxed into the SDP (\ref{eq12}):
\begin{equation}\nonumber
    \min  \langle \bR,\bV\bV^T\rangle
   \text{ s.t }
    \left\{\begin{aligned}
    & \langle F_i, \bV \bV^T \rangle  = 2 - d_i, i = 1,\ldots,n \\
    & \text{diag}(\bV \bV^T) = \onevec_{d+1}
    \end{aligned}\right.
    \label{eq12}
\end{equation}
with $\bR = \begin{pmatrix}
                            \bTheta & \frac{1}{2}\theta\\  
                            \frac{1}{2}\theta^T & 0
                            \end{pmatrix}$.
The diagonal constraint can be rewritten $||V_i|| = 1,  \forall i$. We note $s_1= 1, s_i = s_{i-1}+d_{i-1}, i=2,\ldots,n$. For a given GM variable $x_k$, $s_k$ represents the position of the first row representing $x_k$ in $b, c$ and $e$. The exactly-one constraint $\langle F_i \bV \bV^T \rangle = 2 - d_i$ can be rewritten $V_{d+1}^T (\sum_{j = s_i} ^{s_{i+1}-1} V_j) = 2 - d_i$. 

Instead of doing row-by-row updates, we exploit the property that the matrix $\bTheta$ has zero blocks on the diagonal. For a given GM variable $x_k$, we simultaneously optimize all the rows of $\bV$ that correspond to $x_k$ while keeping all other rows fixed. To simultaneously update all these rows, noting $g_i = \sum_{j=1}^{d+1} R_{i,j}V_j$ as before,  we have to solve:
\begin{equation}
    \min\!\!\! \sum_{i=s_k}^{s_k+d_k-1}\!\! V_i^T g_i \;\text{ s.t.}    
    \left\{\begin{aligned}
    &V_{d+1}^T ( \sum_{j = s_k}^{s_k+d_k-1}V_j)  = 2 - d_k\\
    &||V_i|| = 1, \; s_k\leq i < s_{k+1}
    \end{aligned}\right.
    \label{eq13}
\end{equation}
Notice that since the corresponding $d_k \times d_k$ diagonal block matrix of $\bR$ is all zero, the $g_i$ do not depend on the optimized $V_i$ and can be considered as fixed. In the following, we make the assumption that the vectors $g_i$ and $V_{d+1}$ are never colinear (see~\citep{wang2017mixing}).

By the second constraint, every solution row vector $V_i$ must lie in the unit spherical manifold. Furthermore we observe that:
\begin{theorem}\label{thm1}
At the optimum of (\ref{eq13}), every optimized vector $V_i$ lies in the dimension 2 subspace generated by $g_i$ and $V_{d+1}$.
\end{theorem}

The full proof is provided in Appendix~\ref{app1}. Given that $V_i$ lies on the unit sphere and in this dimension 2 subspace, it is entirely defined by the angle it makes with $V_{d+1}$. With vectors on the sphere, the problem above can be rewritten using trigonometric functions, omitting the $d_i$ norm-1 constraints which are implicitly satisfied. To deal with the remaining exactly-one constraint, we write the Lagrangian dual using multiplier $\lambda$.
\begin{theorem}\label{thm2}
If we denote as $\phi_i$ the angle from $g_i$ to $V_{d+1}$, the dual Lagrangian of problem (\ref{eq13}) is:
\begin{equation*}
     h(\lambda) = -\!\!\!\!\!\!\!\!\sum_{i = s_k}^{s_k+d_k-1} \sqrt{||g_i||^2 + 2\lambda ||g_i|| \cos(\phi_i) + \lambda^2} + (d_k - 2)\lambda
\end{equation*}
\end{theorem}

The full proof is provided in Appendix~\ref{app2}. This is a function on a single dimension, which we can optimize using the Newton method. We introduce a notation for the Newton method: we denote by $f \in \mathcal{C}^{k,p}_L(Q)$ the set of functions $f$ such that:
\begin{itemize}
    \item $f$ is $k$ times continuously differentiable on $Q$.
    \item Its $p^{th}$ derivative is Lipschitz continuous on Q with constant L.
\end{itemize}

\begin{theorem}Rate of convergence of the Newton method \citep{nesterov2018lectures}. Let us consider the problem 
\begin{center}
    $\underset{x \in \mathbb{R}^n}{min} \; f(x)$
\end{center}
Let the function $f(.)$ satisfy the following assumptions:
\begin{itemize}
    \item $f \in \mathcal{C}^{2,2}_L(\mathbb{R}^n)$.
    \item There exists a local minimum of the function $f$ with positive  definite Hessian.
    \item The initial starting point $x_0$ is close enough to $x^*$.
\end{itemize}
Then the Newton's method converges quadratically.
\end{theorem}

To reach quadratic convergence with the Newton method, we need to find a starting point which is close enough to $\lambda^*$. This starting point can be produced by analytically solving the Lagrangian dual of the following relaxation of the BCD problem~(\ref{eq13}):
\begin{equation}
    \min_{}  v^T g \quad\text{s.t}\quad v'^T v  = 2 - d_k,  ||v||^2 = d_k
    \label{eq15}
\end{equation}
where $g^T = [g_1^T \;\ldots\; g_{d_1}^T]$, $v'^T = [V_{d+1}^T  \;\ldots\; V_{d+1}^T]$ and $v^T = [V_1^T \;\ldots\; V_{d_1}^T]$. 

\begin{theorem}\label{thm5}
Let us define $r = 4 - \frac{4}{d_k}$ and $\gamma = ||g||^2 - \frac{1}{d_k}(v'^Tg)^2$. The solution to (\ref{eq15}) is given by $v^* = \alpha v' + \beta g$ with 
\begin{center}
    $\alpha = -\frac{1}{d_k}(\beta v'^Tg + d_k - 2)$ \; $\beta = -\sqrt{\frac{r}{\gamma}}$ 
\end{center}
\end{theorem}

The full proof is provided in Appendix~\ref{app3}.

\begin{theorem}\label{thm6}
The block-coordinate descent with fixed last column $V_{d+1}$ is decreasing.
\end{theorem}

The full proof is provided in Appendix~\ref{app4}. This last theorem ensures that the block-coordinate descent will converge. The convergence speed using random row-by-row updates is linear in the neighborhood of a local optimum when the rank is sufficiently large~\citep{wang2017mixing,erdogdu2018convergence}. 

\paragraph{Computational Complexity.}
The update rule for the Newton method is $\lambda^{j+1} = \lambda^j - [\nabla^2 f(\lambda^j)]^{-1} \nabla f(\lambda^j)$ and the first and second derivative of $f$ are:
    \begin{align}
    f'(\lambda) = & \sum_{i=s_k}^{s_k+d_k-1} \frac{g_i^T V_{d+1}^T + \lambda}{||g_i + \lambda V_{d+1}^T||} - (d_k - 2) \\f''(\lambda) = &\sum_{i=s_k}^{s_k+d_k-1} \frac{||g_i||^2 - (g_i^T V_{d+1}^T)^2}{||g_i + \lambda V_{d+1}^T||^3}
    \end{align}

Given that all vectors here have size $r$, by pre-computing the dot products with $g_i$, computing the first derivative requires $O(d_k r)$ operations. For the second derivative, we can pre-compute the numerators $||g_i||^2 - (g_i^T V_{d+1}^T)^2$ and  evaluating the second derivative is again $O(d_k r)$. Overall, one update of the Newton method will require only $O(d_k r)$ operations, which is also what a single round of row-by-row updates over the $d_k$ rows associated with $x_k$ would require. The BCD updates however benefit from the efficient Newton updates. In our experiments, we observe that only few iterations of the Newton method are needed.

\subsection{Producing an integer primal solution}

To produce a primal integer solution from the optimal solution of the convex relaxation, we use the usual random rounding approach of \citep{goemans1995improved}. We generate a random vector $v_r$ on the unit sphere by sampling each coordinate from a Gaussian distribution and then normalizing. We then assign variable $x_k$ to the value that has the largest scalar product with $v_r$. Following this, the integer solution obtained is submitted to a simple greedy search where, for every variable, the best improving change of value (if any) is applied. This is done repeatedly on every variable until a local minimum is reached. The same process is used for all convex relaxations considered. For every instance we repeat this scheme with 50 random vectors and the best solution is kept.

\section{Experiments}\label{experiments}

All experiments were run on a server equipped with a Xeon\textregistered CPU E5-2680 v3 @ 2.50 GHz and 256GB or RAM running Debian Linux 4.19.98-1. All implementations are sequential and use a single thread of computation.

\subsection{Description of solvers} 

We implemented the row-by-row updates with the dualized exactly-one constraints (LR-LAS) as well as the BCD update method (LR-BCD) in C++ with the Eigen3~\citep{eigenweb} sparse matrix library. We use rank $r=\sqrt{2n}$ by default and tried variants with a low rank (down to $5$).
These implementations are avaible at \url{https://github.com/ValDurante/LR-BCD} and compared with two traditional competitors:  the convergent message passing MAP/MRF algorithm \mbox{TRW-S}, as implemented in Open-GM2~\citep{kappes2015comparative}, an efficient C++ MIT-licensed library, in release 3.3.7, available at \url{https://github.com/opengm}. We used a maximum number of iterations 100000 and a tolerance of $10^{-5}$ with TABLE mode activated. We also compared with the WCSP/CFN algorithm  enforcing Virtual Arc Consistency~\citep{cooper2010soft}, as implemented in toulbar2~\citep{hurley2016multi}, another C++ MIT-licensed library, in release 1.1.1, available at \url{https://github.com/toulbar2/toulbar2}. We used a resolution of $10^{-3}$ (using flags \texttt{-nopre -A -C=1000} to turn off all preprocessing but VAC). A primal solution is obtained by diving with default heuristics and no backtrack, simply by adding flag \texttt{-bt=0} to the previous command line.

Overall, we therefore have four solvers: LR-LAS and LR-BCD (parameterized by their rank), \mbox{TRW-S}, and VAC. We also considered using the interior point method implemented in CVXPY~\citep{agrawal2018rewriting}, but preliminary tests showed that it ran several orders of magnitude slower than the Mixing Method on small problems with 120 binary variables from the BiqMac library at \url{http://biqmac.uni-klu.ac.at/biqmaclib.html}.

\subsection{Description of instances} 

We used hard random instances of MAP/MRF and WCSP/CFN problems. We generated random pairwise problems using the WCSP/CFN generator included in toulbar2 with the \texttt{bin-\{n\}-\{d\}-\{t1\}-\{p2\}-\{seed\}} model. This generates a GM with \texttt{n} variables, \texttt{d} values, \texttt{p2} pairwise functions, each containing a random fraction of potentials equal to 0, the rest being uniformly sampled between 1 and 3. The pairwise functions are sampled uniformly without replacement from the set of all possible edges. We generated different types of random problems, varying along the number of variables (50-100), number of values (3-10) and graph density (sparse-dense), keeping the potential functions in a 50\% tightness regime. Dense problems have a complete graph while sparse problems have the number of pairwise functions set to 4 times the number of variables, giving an average degree of 8.

\subsection{Results}

We present in Table~\ref{table:results} a summary of the performance of all solvers on random instances.
Table~\ref{table:results} shows that approximate LP bounds (VAC/\mbox{TRW-S}) are tighter and faster than SDP bounds on sparse problems while SDP becomes tighter than LP on dense problems. For SDP, LR-BCD is several orders of magnitude faster than LR-LAS, with similar bounds on the dense instances. However, on sparse instances, where \mbox{TRW-S} dominates anyway, LR-BCD lower bounds become looser. This can be explained by the fact that the LR-BCD does not enforce the \textit{gangster} constraint, which is known to tighten SDP relaxations while LR-LAS enforces them.

For upper bounds, none of the four solvers strongly stands out. LR-LAS and LR-BCD provide bounds that are close in general. There are some instances where VAC seems to have better results, especially on sparse instances with large domain sizes while \mbox{TRW-S} is stronger on the dense instances with small domain sizes. However, LR-LAS and LR-BCD SDP-based upper bounds start to dominate on the larger and denser instances.

\setlength{\tabcolsep}{2pt}
\begin{table*}[htbp]
\centering
\begin{tabular}{llSSSSSSSS}
\toprule
\multicolumn{2}{c}{density} & \multicolumn{4}{c}{sparse} & \multicolumn{4}{c}{dense} \\
\cmidrule(r){1-2} \cmidrule(rl){3-6} \cmidrule(l){7-10}
\multicolumn{2}{c}{\texttt{n}}
 &\multicolumn{2}{c}{50} &\multicolumn{2}{c}{100} &\multicolumn{2}{c}{50} &\multicolumn{2}{c}{100} \\
 \cmidrule(r){1-2} \cmidrule(rl){3-4} \cmidrule(rl){5-6} \cmidrule(rl){7-8} \cmidrule(l){9-10}
\multicolumn{2}{c}{\texttt{d}}
 & 3 & 10 & 3 & 10 & 3 & 10 & 3 & 10  \\
\midrule
\multirow{4}{*}{lb} 
 & VAC & 129 &  26 &  255 &  53 &  430 &  37 &  1593 &  77 \\
 & TRW & 129 &  26 &  257 &  54 &  442 &  41 &  1631 &  88 \\
 & BCD & 33 &  -712 &  30 &  -1585 &  953 &  178 &  4260 &  1892 \\
 & LAS & 144 &  -132 &  283 & -172  &  1038 &  274 &  4499 & 2437 \\
\midrule
 \multirow{4}{*}{ub} 
 & VAC & 163 &  77 &  332 &  160 &  1212 &  919 &  5128 &  4200 \\
 & TRW & 170 &  92 &  352 &  177 &  1171 &  894 &  4929 &  8926 \\
 & BCD & 185 & 108  &  369 &  230 &  1153 &  917 &  4922 &  4002 \\
 & LAS & 185 &  108 &  370 & 230 & 1162 &  923 &  4911 & 3999 \\
\midrule
 \multirow{4}{*}{\parbox[t]{5mm}{gap\\(\%)}} 
 & VAC & 21.4 &  67.3 &  23.1 &  67.4 &  64.6 &  96.0 & 68.9 &  97.2 \\
 & TRW & 24.4 &  71.7 &  27.0 &  70.1 &  62.3 &  95.5 &  66.9 &  97.8 \\
 & BCD & 82.7 &  759.5 &  91.9 &  789.8 & 17.4 & 80.4 &  13.4 &  52.7 \\
 & LAS & 22.8 &  221.9 &  23.5 & 175.0 &  10.7 &  70.4 & 8.4 & 39.1  \\
\midrule
  \multirow{4}{*}{\parbox[t]{5mm}{cpu\\(s)}} 
 & VAC & 0.004 &  0.041 &  0.0075 &  0.11 &  0.0098 &  0.40 &  0.034 &  2.6 \\
 & TRW & 0.004 &  0.018 &  0.010 &  0.061 &  0.0048 &  0.35 &  0.025 &  3.3 \\
 & BCD & 0.037 &  0.34 &  0.17 &  2.0 &  0.067 &  0.55 &  0.41 &  3.9 \\
 & LAS & 0.54 &  16 &  5.6 & 67 &  1.2 &  36 &  19 & 219 \\
\bottomrule
\end{tabular}
\captionof{table}{Solver performances: lower bound (lb), upper bound (ub), optimality gap (gap) and cpu-time for all four methods (rows) on all types of random instances (columns).\label{table:results}}
\end{table*}

These results suggest that the best result could be obtained using both SDP and propagation based bounds with just two threads of computations. In this scenario, we would use one thread to execute VAC or \mbox{TRW-S}, and another one to enforce LR-BCD. This combination will provide tighter gaps with enhanced upper and lower bounds along a large spectrum of problems at limited extra computational cost: on the hardest dense problems, where LR-BCD brings the strongest improvments, it is almost as fast as VAC or \mbox{TRW-S}. The use  of LR-LAS to obtain stronger lower bounds is far more demanding in terms of cpu-time: on the hardest dense problems, where it could bring the most significant improvements, it runs 70 times slower than LR-BCD, as well as VAC/\mbox{TRW-S}.

\subsubsection{Constant rank relaxations.}

In Table~\ref{table:lowrank}, we empirically explore the world of very low rank relaxations: the rank is reduced from its theoretically safe value to reach values where the optimum value loses its status of lower bound. Indeed, it is important to recall that for a rank $r \leqslant \sqrt{2n}$ the mixing method is not guaranteed to converge to an optimum of the SDP relaxation. 

Using such low ranks on the rd100-10-dense-0 ($\lceil \sqrt{2n} \rceil = 45$) instance, we see that the LR-BCD method speeds-up as the rank decreases, the value of the primal and integer solutions remaining stable until very low ranks are used: the quality and speed of the method then worsen as the required number of iterations grows with very tiny ranks. LR-LAS also benefits from rank reduction since the value of the integer solution remains stable even with very low ranks. As upper-bounds seems mostly insensitive to rank reduction, rank reduction makes sense if quick reliable upper bounds are sufficient. To benefit from an optimality gap in this setting of very low ranks, a \emph{dual} lower bound would need to be computed.  The acceleration that BCD updates are clearly visible on this instance: $5040$ iterations of LR-LAS corresponds to $1001 \times 5040 = 5045040$ row updates. This is two order of magnitude more than Newton iterations for LR-BCD.

\setlength{\tabcolsep}{2pt}
\begin{table*}[t]
\setlength{\tabcolsep}{6pt}
\centering
\begin{tabular}{lccccccccc}
  \toprule
   & \multicolumn{5}{c}{LR-BCD} & \multicolumn{4}{c}{LR-LAS} \\
  \cmidrule(r){1-1} \cmidrule(rl){2-6} \cmidrule(l){7-10}
 
  \multirow{2}{*}{Rank} & \multirow{2}{*}{value} & \multirow{2}{*}{ub} &  \multirow{2}{*}{cpu (s)} & \multicolumn{2}{c}{iterations} & \multirow{2}{*}{value} & \multirow{2}{*}{ub} & \multirow{2}{*}{cpu (s)} & \multirow{2}{*}{iterations}\\
       &       &       &         & BCD & Newton                  &       &       &         &\\
         \cmidrule(r){1-1} \cmidrule(rl){2-6} \cmidrule(l){7-10}

  45 & 1890 & 4006 & 2.54 & 56 & 12854 & 2438 & 4006 & 212.33 & 5040\\
  41 & 1890 & 4006 & 2.32 & 56 & 12742 & 2440 & 4006 & 192.57 & 4990\\
  37 & 1890 & 4006 & 2.16 & 57 & 13020 & 2443 & 4006 & 181.17 & 5057 \\
  33 & 1890 & 4006 & 1.88 & 55 & 12618 & 2427 & 4006 & 159.39 & 4946\\
  29 & 1890 & 4006 & 1.65 & 54 & 12305 & 2425 & 4006 & 146.72 & 5028\\
  25 & 1890 & 4006 & 1.51 & 56 & 12881 & 2443 & 3972 & 129.25 & 5060\\
  21 & 1890 & 4006 & 1.31 & 56 & 12903 & 2439 & 4006 & 113.51 & 5082\\
  17 & 1890 & 3972 & 1.20 & 61 & 14151 & 2449 & 4006 & 96.12 & 5092\\
  13 & 1890 & 4006 & 1.03 & 64 & 15086 & 2460 & 3972 & 79.59 & 5124 \\
  9 & 1891 & 3972 & 0.93 & 74 & 17795 & 2461 & 4003 & 62.78 & 5116\\
  5 & 1910 & 4006 & 1.26 & 138 & 38037 & 2511 & 4006 & 45.08 & 5148\\
  \bottomrule
\end{tabular}
\captionof{table}{Rank results on the rd100-10-dense-0 instance: optimum value (value), integer solution value (ub), cpu-time (cpu) and number of iterations. For LR-BCD, we report the number of matrix (BCD) and Newton updates. For LR-LAS, we report the number of matrix updates.\label{table:lowrank}}
\end{table*}

\section{Conclusion}

An increasing number of proposals for solving convex relaxations of discrete optimization problems rely on the Burer-Monteiro approach~\citep{burer2003nonlinear}. These approaches are fast and offer both a lower and upper bound which make them highly desirable when efficiency is required. They can also offer tighter gaps than the more usual LP-related local consistencies, especially on the hardest dense problems. However, to the best of our knowledge, all these proposals deal only with the simple case where only diagonal-one constraints exist. This is enough for MaxCUT, MAX2SAT, and MAP/Ising and even MAP/Potts (where potential functions are weighted identity matrices) or even MAXSAT (at the cost of a significantly weakened relaxation for $k$-clauses, $k>2$). 

But a variety of real world problems expressed as optimization over graphical models, in resource allocation, computational biology or image processing, naturally involve variables with many values and arbitrary pairwise functions. This work shows that the extension of the Burer-Monteiro approach to arbitrary domain sizes and functions is not straightforward: even if it improves over usual interior points methods, as implemented in CVXPY, the direct dualization of the required ``exactly-one'' constraints with adequate penalization leads to slow convergence, removing most of the practical interest of the Burer-Monteiro approach. Dedicated ways of dealing with these exactly-one constraints are needed. The BCD approach introduced is a clear step forward in this direction. With its fast Newton updates,  it offers important speed-ups compared to the dualized variant, at the sole cost of relaxing the gangster constraint. On the hardest instances, where VAC/\mbox{TRW-S} can only provide very large gaps, the BCD approach we described is several orders of magnitude faster than interior points methods and orders of magnitude faster than the Lassere-dualized approach, while offering much tighter optimality gaps than VAC/\mbox{TRW-S}. It becomes the method of choice for hard problems and can be used in parallel with the linear bounds produced by VAC/\mbox{TRW-S} to get the best of both worlds, at modest additional runtime. 

The possibility of relying on extremely low rank is another attractive capacity of the Burer-Monteiro method. We observe that the CPU-time that is required to produce an upper bound can be easily further reduced. A dual bound would then be able to produce an optimality gap that could exploited in the context of Branch-and-Bound methods, providing a desirable tight bounding method with an adjustable cpu/tightness compromise.

\section{Acknowledgments and Disclosure of Funding}
This work has benefited from the AI Interdisciplinary Institute
ANITI funding, through the French “Investing for the Future PIA3” program under the Grant agreement n°ANR-19-PI3A-0004. The first author is also supported by the EUR BioEco program. The second and third authors are both supported by the "Décomposition de Modèles Graphiques – DE-MO-GRAPH" program under grant n°ANR-16-CE40-0028.

\printbibliography

\newpage

\appendix

\section{Appendix}

\subsection{Proof of Theorem~\ref{thm1}}\label{app1}

\begin{proof}
Let $V_i^*$ be the value of $V_i$ in an optimum of the problem above and assume that $V_i^*$ does not lie in the dimension 2 subspace generated by $g_i$ and $V_{d+1}$. Let $u=V_{d+1}, v$ be the orthonormal basis of this subspace such that $g_i$ has a positive coordinate over $v$. Then, $V_i$ can be written as $V_i= \alpha u + \beta v + \gamma w$ where $w$ lies in the $r-2$ dimensional subspace orthogonal to the subspace generated by $g_i$ and $V_{d+1}$ and $\gamma \neq 0$. Since, when $\beta>0$,  changing the sign of $\beta$ improves the criteria without changing the status of the ``exactly-one'' and norm-1 constraints, we must have $\beta<0$. Let $\beta' = -\sqrt{\beta^2+\gamma^2}$ and $V_i^*{}' = \alpha u + \beta' v$. If we now consider a solution where $V_i^*$ has been replaced by $V_i^*{}'$, the ``exactly one'' constraint is still satisfied because the scalar product of $V_i^*{}'$ with $V_{d+1}$ is unchanged. The norm $1$ constraint is still satisfied as $||V_i^*{}'||^2= \alpha^2+\beta'^2 = \alpha^2+\beta^2+\gamma^2= 1$. Then, the new solution improves the criteria given that $\beta' < \beta < 0$ and that $g_i^T v > 0$.
\end{proof}

\subsection{Proof of Theorem~\ref{thm2}}\label{app2}

\begin{proof}
With $V_i$ being norm one and lying in the dimension 2 subspace generated by $g_i$ and $V_{d+1}$, let $\psi_i$ be the angle from $V_{d+1}$ to $V_i$ in this subspace. We can rewrite the problem (\ref{eq13}) as:

\begin{equation}
    \min_{\psi_i\in[0,\pi]} \sum_{i=s_k}^{s_k+d_k-1}||g_i||\cos(\phi_i+\psi_i) \quad\text{s.t} \quad    
    \sum_{i=s_k}^{s_k+d_k-1} \cos(\psi_i)  = 2 - d_k
\end{equation}
The Lagrangian is: 
$$L(\lambda,\Psi) = \sum_{i=1}^{d_1}||g_i||\cos(\phi_i+\psi_i) + \lambda(\sum_{i=s_k}^{s_k+d_k-1} \cos(\psi_i) + d_k - 2)$$

The partial derivatives of $L$ are $\pdv{L}{\psi_i} = -||g_i||\sin(\phi_i+\psi_i) - \lambda \sin(\phi_i)$ which must all be equal to zero at the optimum. We use the fact that the sum of sinusoids with the same period and different phases is also a sinusoid:

$$ A\sin(\omega t+ \phi) + B \sin(\omega t) = \sqrt{A^2+B^2+2AB\cos(\phi)} \sin(\omega t + \arctan(\frac{A\sin(\phi)}{A\cos(\phi)+B})) $$
We have:
$$-||g_i||\sin(\phi_i+\psi_i) - \lambda \sin(\phi_i) = 
   \sqrt{||g_i||^2+\lambda^2+2\lambda||g_i||\cos(\phi_i)}\sin(\psi_i+\arctan(\frac{||g_i||\sin(\phi_i)}{||g_i||\cos(\phi_i)+\lambda}))$$
which implies that
$$ \psi_i + \arctan(\frac{||g_i||\sin(\phi_i)}{||g_i||\cos(\phi_i)+\lambda}) = \pm \pi $$

Since $-\frac{3}{2}\pi < -\pi - \arctan(\frac{||g_i||\sin(\phi_i)}{||g_i||\cos(\phi_i)+\lambda}) < -\frac{\pi}{2} < 0$, we have:
$$ \psi_i =\pi - \arctan(\frac{||g_i||\sin(\phi_i)}{||g_i||\cos(\phi_i)+\lambda}) $$
By trigonometric identities, noting $\gamma_i = \frac{||g_i||\sin(\phi_i)}{||g_i||\cos(\phi_i)+\lambda}$, we have:
$$ \cos(\psi_i) = -\cos(\arctan(\gamma_i)) = -\frac{1}{\sqrt{1+\gamma_i^2}}$$

Developing and simplifying $(1+\gamma_i^2)$ and plugging the result back above, we get:
$$ \cos(\psi_i) = -\frac{||g_i||\cos(\phi_i)+\lambda}{||g_i +\lambda V_{d+1}||}$$

Similarly, one can derive:
$$\cos(\phi_i+\psi_i) = -\frac{||g_i||+\lambda\cos(\phi_i)}{||g_i+\lambda V_{d+1}||}$$

Plugging this back into the Lagrangian, we get:

\begin{eqnarray*}
L(\lambda,\Psi) & = & \sum_{i=s_k}^{s_k+d_k-1} \left(-||g_i||\frac{||g_i|| + \lambda\cos(\phi_i)}{||g_i+\lambda V_{d+1}||}\right) + 
                     \lambda\left(\sum_{i=s_k}^{s_k+d_k-1} \left( -\frac{||g_i||\cos(\phi_i)+\lambda}{||g_i +\lambda V_{d+1}||}\right) + d_k - 2\right)\\
                & = & -\sum_{i=s_k}^{s_k+d_k-1}(||g_i+\lambda V_{d+1}||) + \lambda (d_k-2)
\end{eqnarray*}

So the dual problem is to find $\lambda$ that maximizes $h(\lambda) = -\sum_{i=s_k}^{s_k+d_k-1}(||g_i+\lambda V_{d+1}||) + \lambda (d_k-2)$. 
The first derivative of $h$ is:
$$h'(\lambda) = -\sum_{i=s_k}^{s_k+d_k-1}\left(\frac{g_i^TV_{d+1}+\lambda)}{||g_i+\lambda V_{d+1}||}\right) + (d_k-2)$$
and the second derivative:
$$h''(\lambda) = -\sum_{i=s_k}^{s_k+d_k-1}\left(\frac{||g_i||^2-(g_i^TV_{d+1})^2}{||g_i+\lambda V_{d+1}||^3} \right)$$

One can see that $h''(\lambda) < 0$ whenever $g_i$ and $V_{d+1}$ are not colinear which proves $h$ is strictly concave and guarantees the unicity of the optimum. Being strictly concave, we just have to find $\lambda^*$ such that $h'(\lambda^*) = 0$.
\end{proof}

\subsection{Proof of Theorem~\ref{thm5}}\label{app3}

Let us first introduce a lemma that will be useful for the proof of this theorem.
\begin{lemma}
The solution $v^*$ to (\ref{eq15}) must lie in the subspace $vec(v',g)$.
\end{lemma}

Let us consider $P_1 : \big \{v \in \mathbb{R}^n \; | \; v'^Tv = 2 - d_k  \big \}$ and $D_1 : \big \{v \in \mathbb{R}^n \; | \; ||v||^2 = d_k  \big \}$. For $y \in P_1$, we have $y = y_0 + z$ with $y_0$ a particular solution to the linear equation $v'^Ty_0 = 2 - d_k$ and $z$ a vector such that $v'^Tz = 0$. For $y_0$, we take $y_0 = \frac{2-d_i}{||v'||^2}v'$ which gives us the desired solution $v'^Ty_0 = \frac{2-d_k}{||u||^2}v'^Tv' = 2- d_k$. Since $y \in D_1$, the squared norm of $z$ must also satisfies $||z||^2 = d_k - (2-d_k)^2$.\\
Then, we compute the dot product $g^Ty = \frac{(2-d_k)}{||v'||^2}v'^Tg + z^Tg$. The first term is constant for all $y$ that lies in the feasible set $P_1 \cap D_1$. We just have to check the value of the second term $z^Tg$. Let $\Pi_{P_1}(g)$ be the projection of $g$ into the hyperspace orthogonal to $v'$. We can write the dot product : $z^Tg = z^T(g - \Pi_{P_1}(g) + \Pi_{P_1}(g)) = z^T(g - \Pi_{P_1}(g)) + z^T\Pi_{P_1}(g) = z^T\Pi_{P_1}(g)$ since $z^T(g - \Pi_{P_1}(g)) = 0$. Thus we only worry about the value of the dot product $z^T\Pi_{P_1}(g)$ which is minimal in the direction $-\Pi_{P_1}(g)$. Since $\Pi_{P_1}(g)$ lies in the space generated by $v'$ and $g$ we proved that the solution to (\ref{eq15}) must lie in the space generated by $v'$ and $g$.

\begin{proof}
Let us consider $v \in vec(v',g)$ e.g. $v = \alpha v' + \beta g$ for some $\alpha, \beta \in \mathbb{R}$. The first constraint gives us:
\begin{center}
    $\alpha = -\frac{1}{d_k}(\beta v'^Tg + d_k - 2)$
\end{center}
By developing the squared norm we have:
\begin{center}
    $\beta^2(\frac{1}{d_k}(v'^Tg)^2 - \frac{2}{d_k}(v'^Tg)^2 + ||g||^2) + \frac{(d_k-2)^2}{d_k} = d_k $\\
    \vspace{0.25cm}
    $\beta^2(||g||^2 - \frac{1}{d_k}(u^Tg)^2) = d_k - \frac{(d_k-2)^2}{d_k}$\\
    \vspace{0.25cm}
    $\beta^2(||g||^2 - \frac{1}{d_k}(v'^Tg)^2) = 4 - \frac{4}{d_k}$
\end{center}
Using the Cauchy-Schwartz inequality we know that $||g||^2 - \frac{1}{d_k}(v'^Tg)^2 > 0$ whenever $v'$ and $g$ are not colinear. Thus using the notation $r = 4 - \frac{4}{d_k}$:
\begin{center}
$\beta = \mp \sqrt{\frac{r}{\gamma}}$ with $\gamma = ||g||^2 - \frac{1}{d_k}(v'^Tg)^2$    
\end{center}
One can easily see that the objective value is smaller for $\beta = -\sqrt{\frac{r}{\gamma}}$ which gives us the result.
\end{proof}

\subsection{Proof of Theorem~\ref{thm6}}\label{app4}

\begin{proof}  Let us show that the objective difference before and after a BCD update of $V$ is always positive. For a given column $V_i$ the objective difference before and after the update is:
\begin{center}
    $f(V_i) - f(\hat{V_i}) = g_i^T(V_i - \hat{V_i})$\\
    \vspace{0.25cm}
    with $\hat{V_i} = \alpha_i V_{d+1} + \beta_i g_i$\\
    \vspace{0.25cm}
    $\Rightarrow g_i = \frac{1}{\beta_i}(\hat{V_i} - \alpha_i V_{d+1})$\\
    \vspace{0.25cm}
    $g_i^T(V_i - \hat{V_i}) = \frac{1}{\beta_i}(\hat{V_i} - \alpha_i V_{d+1})^T(V_i - \hat{V_i})$\\
    \vspace{0.25cm}
    $= \frac{1}{\beta_i}\big( (\hat{V_i}^TV_i - 1) - \alpha_i V_{d+1}^T(V_i - \hat{V_i})\big) $
\end{center}
Let us first check the first term. We know that:
\begin{center}
    $\frac{1}{\beta_i}(\hat{V_i}V_i - 1) = -\frac{1}{2\beta_i}||V_i - \hat{V_i}||^2$
\end{center}
Moreover: 
\begin{center}
    $\beta_i = \frac{||g_i|| cos(\phi_i + \psi_i) - cos(\psi_i)V_{d+1}^Tg_i}{||g_i||^2 - (V_{d+1}^Tg_i)^2}$
\end{center}
Using the Cauchy-Schwartz inequality, we already know that $||g_i||^2 - (V_{d+1}^Tg_i)^2 \geqslant 0$.
Next we will rewrite the numerator:
\begin{center}
$||g_i|| cos(\phi_i + \psi_i) - cos(\psi_i)V_{d+1}^Tg_i = ||g_i||(cos(\phi_i + \psi_i) - cos(\psi_i)cos(\phi_i))$\\
\vspace{0.25cm}
$ = -||g_i||sin(\phi_i)sin(\psi_i)$
\end{center}
Since $\phi_i \in [0,\pi]$ and $\psi_i \in [0,\pi]$ the sinus are positive so we can conclude that $\beta_i \leqslant 0$. We proved that the first term is positive, let us now check the second term.
\begin{center}
    $-\alpha_i V_{d+1}^T(V_i - \hat{V_i}) = - \alpha_i V_{d+1}^T(V_i - (\alpha_i V_{d+1} + \beta_i g_i))$\\
    \vspace{0.25cm}
    $= -\alpha_i V_{d+1}^TV_i + \alpha_i^2 + \alpha_i \beta_i V_{d+1}^Tg_i$.
\end{center}
Using $\alpha_i = V_{d+1}^TV_i - \beta_i V_{d+1}^Tg_i$ we have:
\begin{center}
    $-\alpha_i V_{d+1}^T(V_i - \hat{V_i}) = - (V_{d+1}^TV_i - \beta_i V_{d+1}^Tg_i)V_{d+1}^TV_i + \alpha_i^2 + (V_{d+1}^TV_i - \beta_i V_{d+1}^Tg_i)\beta_i V_{d+1}^Tg_i$\\
    \vspace{0.25cm}
    $= -(V_{d+1}^TV_i)^2 + \beta_i(V_{d+1}^Tg_i)(V_{d+1}^TV_i) + \alpha_i^2 + \beta_i(V_{d+1}^TV_i)(V_{d+1}^Tg_i) - \beta_i^2(V_{d+1}^Tg_i)^2$\\
    \vspace{0.25cm}
    $ = 2\beta_i(V_{d+1}^TV_i)(V_{d+1}^Tg_i) - \beta_i^2(V_{d+1}^Tg_i)^2 - (V_{d+1}^TV_i)^2 + \alpha_i^2$\\
    \vspace{0.25cm}
    $= - (\beta_i(V_{d+1}^Tg_i) - V_{d+1}^TV_i)^2 + (V_{d+1}^TV_i - \beta_i(V_{d+1}^Tg_i))^2 = 0$
\end{center}
Thus, the second term is equal to 0 so we proved that $f(V_i) - f(\hat{V_i}) \geqslant 0$
\end{proof}

\end{document}